\documentclass[11pt,
]{article}

\usepackage{amssymb,amstext,amsthm,latexsym}
\usepackage{amsmath}
\usepackage{amsfonts}
\usepackage{mathrsfs}
\usepackage{mparhack}
\usepackage[bbgreekl]{mathbbol}

\newcommand{\nek}{\newcommand}
\nek{\renek}{\renewcommand}

\makeatletter
\if@twoside%
\else\fi%
\makeatother

\nek{\punk}{\subsection}

\nek{\skr}{\mathscr}
\nek{\cP}{\skr P}

\nek{\vyk} [1] {}

\nek{\imar}[1] 
{\marginpar[
\flushright\footnotesize%
$\mtho\longrightarrow$\\%
\vspace{-1ex}{#1}\vspace*{1ex}]%
{
\flushleft\footnotesize%
$\mtho\longleftarrow$\\%
\vspace{-1ex}{#1}\vspace*{1ex}}%
}%

\renek{\imar}[1]{}
\nek{\imae} {\imar}

\nek{\itsep}{\itemsep=0.3ex plus 0.1ex minus 0.1ex}
\nek{\tenu}[1]{
\def\theenumi{#1}
\itsep
}

\theoremstyle{plain}

\newtheorem{theore}             {Theorem} [subsection] 
\newtheorem{corollar}  [theore]{Corollary}
\newtheorem{propo}     [theore]{Proposition}
\newtheorem{lemm}      [theore]{Lemma}
\newtheorem{lemt}      {Lemma}[theore]
\newtheorem{cla}       [theore]{Claim} 
\newtheorem{clt}       {Claim} [theore]

\theoremstyle{definition}
\newtheorem{defn}      [theore]{Definition}
\newtheorem{rem} [theore]   {Remark} 
\newtheorem*{solu}{Solution}      
\newtheorem*{probl}{Problem}      
\newtheorem{exz} {Example} 
\newtheorem*{prF}{Proof}               

\nek{\thsp}{\hspace{0.1ex plus \mathsurround}}

\nek{\bex}{\begin{exz}}
\nek{\eex}{\end{exz}}
\nek{\bprl}{\begin{probl}}
\nek{\eprl}{\end{probl}}
\nek{\bsol}{\begin{solu}}
\nek{\esol}{\end{solu}}
\nek{\bpro}{\begin{propo}}
\nek{\epro}{\end{propo}}
\nek{\bcor}{\begin{corollar}}
\nek{\ecor}{\end{corollar}}
\nek{\bdf} {\begin{defn}} 
\nek{\eDf} {\end{defn}}
\nek{\edf} {\end{defn}}
\nek{\edF} {\end{defn}}
\nek{\ble} {\begin{lemm}}
\nek{\ele} {\end{lemm}}
\nek{\blt} {\begin{lemt}}
\nek{\elt} {\end{lemt}}
\nek{\bre}{\begin{rem}}
\nek{\ere}{\end{rem}}
\nek{\bte} {\begin{theore}}
\nek{\ete} {\end{theore}}
\nek{\bpf} {\begin{prF}} 
\nek{\epf} {\qed\end{prF}} 
\nek{\ePf} {\end{prF}} 
\nek{\qeD} {\qed}
\nek{\qeDD} [1] 
{\hfill\hbox{\qed~({\small #1\/}\hspace{0.1ex})}}
\nek{\epF} [1] {\qeDD{#1}\end{prF}} 

\nek{\ben}{\begin{enumerate}\itsep}
\nek{\een}{\end{enumerate}}
\nek{\bit}{\begin{itemize}\itsep}
\nek{\eit}{\end{itemize}}
\nek{\bde}{\begin{description}\itsep}
\nek{\ede}{\end{description}}
\nek{\bay}{\begin{array}}
\nek{\eay}{\end{array}}
\nek{\bce}{\begin{center}}
\nek{\ece}{\end{center}}
\nek{\bqo}{\begin{quotation}\noi}
\nek{\eqo}{\end{quotation}}
\nek{\ZFC}{{\bf ZFC}}
\nek{\ZF}{{\bf ZF}}

\nek{\iesp}{\hspace{0.3ex}}
\nek{\resp}{\hspace{0.25ex}}
\nek{\ie} {\hbox{\sl i.\iesp e.}}
\nek{\eg} {\hbox{\sl e.\iesp g.}}
\nek{\noo}
{\hbox{\sl w.\iesp l.\iesp o.\iesp g.}}
\nek{\vrt} {\hbox{w.\iesp r.\iesp t.}}
\nek{\ddd}[1]{$\mtho\hspace{0.2ex}{#1}$-\hspace{0.0ex}}
\nek{\dd}{\ddd}
\nek{\ran}  {\mathop{\tt ran}}
\nek{\otp} [1] {\mathop{\rm otp}(#1)}
\nek{\dom}  {\mathop{\tt dom}}
\nek{\tsup} {\mathop{\tt sup}}
\nek{\tmax} {\mathop{\tt max}}
\nek{\tlim} {\mathop{\tt lim\hspace{0.3ex}}}

\nek{\al} {\alpha}
\nek{\La}{\Lambda}
\nek{\la}{\lambda}
\nek{\sg} {\sigma}
\nek{\Sg} {\Sigma}
\nek{\ve}{\varepsilon}
\nek{\vpi}{\varphi}
\nek{\vt} {\vartheta}
\nek{\om} {\omega}
\nek{\Om} {\Omega}
\nek{\lom}{^{<\om}}
\nek{\omi} {\om_1}
\nek{\ali} {\aleph_1}
\nek{\alo} {{\aleph_0}}

\nek{\fs}[2]{{\hspace*{0.3ex}\boldsymbol\Sigma}^{#1}_{#2}}
\nek{\fp}[2]{{\boldsymbol\Pi}^{#1}_{#2}}
\nek{\fd}[2]{{\boldsymbol\Delta}^{#1}_{#2}}

\nek{\id}[2]{{\varDelta}^{#1}_{#2}}
\nek{\ip}[2]{{\varPi}^{#1}_{#2}}
\nek{\is}[2]{{\varSigma}^{#1}_{#2}}
\nek{\iSg}{\varSigma}
\nek{\iPi}{\varPi}

\nek{\BBB}{\hspace{0.01ex}}
\nek{\dR}{{\BBB{\mathbb R}\BBB}}

\nek{\bn}{\om^\om}
\nek{\sus} {\mathopen{\exists\hspace{0.35ex}}}
\nek{\kaz} {\mathopen{\forall\hspace{0.35ex}}}
\nek{\imp} {\Longrightarrow} 
\nek{\mpi} {\Longleftarrow} 
\nek{\eqv} {\Longleftrightarrow} 
\nek{\leqv} {\;\eqv\;} 
\nek{\limp} {\;\imp\;} 
\nek{\ti}  {\times} 
\nek{\sq}  {\subseteq}
\nek{\su}  {\subset}
\nek{\sneq}{\subsetneqq}
\nek{\we}  {{\mathbin{\hspace{0.15ex}^\wedge}}}
\nek{\obr} {^{-1}}
\nek{\dif} {\smallsetminus}
\nek{\res} {\mathbin{\restriction}}
\nek{\lef} {\preccurlyeq}
\nek{\gef} {\succcurlyeq}
\nek{\pu}  {\varnothing}
\nek{\iy}  {\infty}
\nek{\piy} {+\iy}
\nek{\nin} {\not\in}
\nek{\onto}{\stackrel{\text{\rm onto}}{\longrightarrow}}
\nek{\ang} [1] {\langle #1\rangle}
\nek{\stk} [2] {\ang{#1\hspace{0.3ex};\hspace{0.1ex}#2}}
\nek{\sis} [2] {\ans{#1}_{#2}}
\nek{\ans} [1] {\{\hspace{0.01ex}#1\hspace{0.01ex}\}}

\nek{\zz} {\linebreak[0]} 
\nek{\ens} [2] {\ans{{#1\hspace{0.5ex}{:}}\zz\hspace{0.5ex}#2}}
\nek{\itla} {\item\label}

\nek{\ubf}{\fontseries{b}\selectfont}
\nek{\TS}{\textstyle}
\nek{\yo} {,\linebreak[0]}
\nek{\yi} {\hspace{0.2ex},\linebreak[0]\hspace{0.2ex}}
\nek{\yd} {\hspace{0.2ex},\linebreak[0]\:}
\nek{\yt} {\hspace*{0.2ex},\linebreak[0]\;} 
\nek{\lex} {<_{\text{\tt lex}}}
\nek{\lexe} {\leqslant_{\text{\tt lex}}}

\nek{\snos} [1] {\,\footnote{\hspace*{2pt}#1}}
\nek{\renu}{\tenu{{\rm(\roman{enumi})}}}
\nek{\fenu}{\tenu{{\rm(\fnsymbol{enumi})}}}
\nek{\Renu}{\tenu{{\rm(\Roman{enumi})}}}
\nek{\rit} [1] {{\it#1\/}}
\nek{\lam} [1] {\label{#1}\hspace*{-3pt}\imar{#1}}%
\nek{\las} [1] {\label{#1}\imar{#1}}%

\nek{\atc}{\addtocounter{enumi}{1}}
\nek{\wo} {\text{\ubf WO}}
\nek{\bez} {\dif}
\nek{\AC}   {{\text{\bf AC}}}
\nek{\DC}   {{\text{\bf DC}}}
\nek{\ROD}  {{\text{\bf{ROD}}}}
\nek{\hrod}{{\text{\bf HROD}}}
\nek{\OD} {{\text{\bf OD}}}

\nek{\rL} {\text{\ubf L}}
\nek{\rV} {\text{\ubf V}}
\nek{\rU} {\text{\ubf U}}

\nek{\led} {\leq^\ast}
\nek{\ld} {<^\ast}

\nek{\pqo} {PQO}

\nek{\Ord}  {\mathop{\tt Ord}}
\nek{\supp}  {\mathop{\tt supp}}
\nek{\card}  {\mathop{\tt card}}

\nek{\ba} {\beta}
\nek{\ga} {\gamma}
\nek{\ka} {\kappa}

\nek{\np}{\newpage}

\nek{\mtho}{\mathsurround=0mm}
\nek{\msur}{\hspace{-1\mathsurround}}
\nek{\noi}{\noindent}
\nek{\vom}{\vspace{1mm}}
\nek{\vim}{\vspace{-1mm}}
\nek{\vtm}{\vspace{2mm}}

\nek{\eqr} {equivalence relation}

\nek{\rE} {\mathrel{\mathsf E}}

\nek{\ek}[2] {[#1]_{{#2}}}

\nek{\eke}[1] {\ek{#1}{\rE}}

\nek{\qand}{\quad\text{and}\quad}

\nek{\dX}{{\BBB{\mathbb X}\BBB}}
\nek{\dN}{{\BBB{\mathbb N}\BBB}}
\nek{\dB}{{\BBB{\mathbb B}\BBB}}
\nek{\dF}{{\BBB{\mathbb F}\BBB}}
\nek{\fun}{{\BBB{\mathbb {Fun}}\BBB}}

\nek{\doP}  [1] {{#1}^\complement}

\nek{\curle}{\preccurlyeq}
\nek{\cle}{\curle}
\nek{\cl} {\prec}
\nek{\ncle}{\not\cle}
\nek{\ncl}{\not\cl}

\nek{\gh}{\mathbb P}
\nek{\ghd}{\gh^2}

\nek{\fdt}{\hbox{\raisebox{-0.25ex}{\LARGE\bf.}}}
\nek{\bdot}[1] 
{\raisebox{-0.07ex}{\mtho$\stackrel{\fdt}{#1}$}}

\nek{\dox}{\bdot{\text{\bf x}}}
\renek{\dox}{\bdot{\boldsymbol x}}
\nek{\doxl}{\dox_{\tt le}}
\nek{\doxr}{\dox_{\tt ri}}

\nek{\dn}{2^\om}
\nek{\dpe} {\mathord{\drof\gh\rE}}
\nek{\dpew}{\mathord{\drow{\gh}\rE}}

\nek{\trof} [3] {{#1}\ti_{#2}{#3}}

\nek{\drof} [2] {{#1}\ti_{#2}{#1}}
\nek{\drow} [2] {{#1}\ti_{#2}^{\text{weak}}{#1}}

\nek{\ck} {\om_1^{\text{\sc ck}}}

\nek{\Eo}  {\rE_{\text{\sf0}}}
\nek{\Fo}  {\rF_{\text{\sf0}}}
\nek{\nE}  {\mathbin{{\not\hspace{-0.35ex}\sf E}}}

\nek{\bcl} {\begin{cla}}
\nek{\ecl} {\end{cla}}
\nek{\bct} {\begin{clt}}
\nek{\ect} {\end{clt}}

\nek{\bV}{{\mathbf V}}
\nek{\bL}{{\bf L}}

\nek{\gM} {\mathfrak M}
\nek{\cM} {\skr M}

\nek{\ccs} {}
\nek{\cF}{{\ccs{\skr F}\ccs}}
\nek{\cS}{{\ccs{\skr S}\ccs}}
\nek{\cD}{{\ccs{\skr D}\ccs}}

\nek{\cf} [1] {\cF_{#1}}
\nek{\cfx} [1] {\cF_#1}
\nek{\cfd} [2] {\cF_{#1}(#2)}

\nek{\etc} {{\sl etc}}

\nek{\bus}{\begin{equation}}   
\nek{\eus}{\end{equation}}

\nek{\pp} [2] {\gh^{#1}_{#2}}

\nek{\dpd} [1] {\gh\ti_{#1}\gh}
\nek{\dpx} [1] {\gh\ti_{\rE_{#1}}\gh}

\nek{\dpw} {\gh(W)}
\nek{\dpwe} {\gh(W)\ti_{\rE}\gh(W)}

\nek{\ups} {\varUpsilon}

\nek{\apr} {\approx}
\nek{\napr}{\not\apr}

\nek{\gp}{\mathfrak p}
\renek{\gp}{\mathbb p}

\nek{\aenu}{\tenu{{\rm(\arabic{enumi})}}}
\nek{\Aenu}{\tenu{{\rm(\Alph{enumi})}}}
 
\nek{\doxlp}{\dox{}'_{\tt le}}
\nek{\doxrp}{\dox{}'_{\tt ri}}

\nek{\bfit}{\bfseries\itshape}

\nek{\esn}{\rS_{\ans{1/n}}}
\nek{\rS}  {\mathbin{\sf S}}

\nek{\wh}{\widehat}
\nek{\wY}{\widehat Y}
\renek{\wY}{C}

\nek{\lr}{LR}
\nek{\rl}{RL}

\nek{\vx} {\vec x}

\nek{\Xa}{X^*}
\nek{\Ua}{U^*}
\nek{\Va}{V^*}
\nek{\Do}{D_0}

\nek{\rEF} {\mathrel{\mathsf E_\cF}}
\nek{\refx} [1] {\mathrel{\mathsf E_{\cfx#1}}}

\begin{document}

\title
{Bounding and decomposing thin analytic partial orderings
\thanks{Partial support of ESI 2013 Set theory program, Caltech, and 
RFFI grant 13-01-00006 acknowledged.}}

\author 
{Vladimir~Kanovei\thanks{IITP RAS and MIIT,
  Moscow, Russia, \ {\tt kanovei@googlemail.com}}  
}

\date 
{\today}

\maketitle

\begin{abstract}
We modify arguments in \cite{hms} to reprove extensions of 
two key results there in the context of bounding and decomposing 
of analytic subsets of 
Borel partial quasi-orderings.
\end{abstract}

\subsection{Introduction}
 
The following theorem is the main content of this note.

\bte
\lam{mt}
Let\/ $\cle$ be a\/ $\id11$ PQO on\/ $\bn$, 
$\apr$ be the associated \eqr, 
and\/ $\Xa\sq \bn$ be a\/ $\is11$ set such that\/ 
${\cle}\res \Xa$ is thin\snos
{Meaning that there is no perfect set of pairwise 
\dd\cle incomparable elements.}.
Then\/ 
\ben
\renu
\itla{mtb} 
there is an ordinal\/ $\al<\ck$ and a\/ $\id11$ \lr\ order 
preserving map\/ $F:\stk{\bn}{\cle}\to\stk{2^\al}{\lexe}$ 
satisfying the following additional requirement$:$ 
if\/ $x,y\in\Xa$ then\/ ${x\napr y}\limp{F(x)\ne F(y)}\;;$ 

\itla{mtd} 
$\Xa$ is covered by the countable union of\/ 
all\/ $\id11$ \dd\cle chains\/ $C\sq\bn\;.$ 
\een
\ete 

The theorem is essentially proved in \cite[3.1 and 5.1]{hms}. 
Literally, only the case of $\id11$ subsets $\Xa$ is considered 
in \cite{hms}, but the case of $\is11$ sets $\Xa$ can be 
obtained by a rather transparent rearrangement of the arguments 
in \cite{hms}.
See also \cite{k:blin} in matters of the additional requirement in 
claim \ref{mtb} of the theorem, 
which also is presented in \cite{hms} implicitly.
Our proofs will largely follow the arguments in \cite{hms}, 
but by necessity we modify those here and there in order to 
streamline some key arguments.
On the other hand, we substitute reflection arguments in \cite{hms} 
with more transparent constructions.

\subsection{Notation and an important lemma}
\las{nota}

\paragraph{\normalsize\ubf Non-strict relations}
\bde
\item[\rm PQO, \it partial quasi-order\/$:$] \ 
reflexive ($x\le x$) and transitive in the domain;


\item[\rm LQO,  
\it linear quasi-order\,$:$] \  
PQO and $x\le y\lor y\le x$  
in the domain; 

\item[\rm LO, \it linear order\,$:$]  \ 
LQO and $x\le y\land y\le x\imp x=y$; 

\item[\it associated equivalence relation\,$:$]  \ 
$x\approx y$ iff $x\le y\land y\le x$.
\ede

\paragraph{\ubf Strict relations}
\bde
\item[\it strict \rm PQO\,$:$] \ irreflexive ($x\not<x$) and
transitive; 

\item[\it strict \rm LQO\,$:$] \  
strict PQO and $x<y\imp\kaz z\:(z<y\lor x<z)$; 

\item[\it strict \rm LO\,$:$]  \ 
strict PQO and the trichotomy $\kaz x,y\:(x<y\lor y<x\lor x=y)$.
\ede
By default we consider only \rit{non-strict} orderings.
All cases of consideration of \rit{strict} PQOs
will be explicitly specified.

Any non-strict PQO $\le$ defines an associated strict
one so that $x<y$ iff $x\le y\land{y\not\le x}$.
In the opposite direction, given a strict PQO $<$, we 
define an \eqr\ $x\approx y$ iff 
$x<z\eqv y<z$ and $z<x\eqv z<y$ for all $z$ in the domain, 
and then define $x\le y$ iff $x<y$ or
$x\approx y$.

\paragraph{\ubf Order preserving maps}
\bde
\item[\it \lr\ (left--right) order preserving map$:$]
any map $f:\stk{X}{\le}\to\stk{X'}{\le'}$ such that 
we have $x\le y\imp f(x)\le' f(y)$ for all $x,y\in\dom f$;

\item[\it \rl\ (right--left) order preserving map$:$]
any map $f:\stk{X}{\le}\to\stk{X'}{\le'}$ such that 
we have $x\le y\mpi f(x)\le' f(y)$ for all $x,y\in\dom f$;

\item[\it 2-ways order preserving map$:$]
any map $f:\stk{X}{\le}\to\stk{X'}{\le'}$ such that 
we have $x\le y\eqv f(x)\le' f(y)$ for all $x,y\in\dom f$.
\ede

\paragraph{\ubf Varia}
\bde
\item[\it sub-order\,$:$] \ 
restriction of the given PQO to a subset of
its domain.

\item[\it $\lex\yt\lexe\;:$]
the lexicographical LOs on sets of the form $2^\al\yd\al\in\Ord$, 
resp.\ strict and non-strict;

\item[\it a\/ \dd\le chain in a PQO $:$]
any set of 2wise \dd\le comparable elements, \ie, LQO;

\item[\it a\/ \dd\le thin set in a PQO$:$]
any set in the domain of $\le$ 
containing no perfect subsets of 2wise \dd\le incomparable elements;


\item[\rm $\eke x=\ens{y\in\dom\rE}{x\rE y}$
(the \dd\rE\rit{class} of $x$) 
and
$\eke X=\textstyle\bigcup_{x\in X}\eke x$]
--- whenever  $\rE$ is an \eqr\ and $x\in\dom\rE$, $X\sq\dom\rE$.
\ede

\ble
[Kreisel selection]
\lam{ks}
Let\/ $D$ be the 
set of all\/ $\id11$ points in\/ $\bn.$
If\/ $P\sq\bn\ti D$ is a\/ $\ip11$ set, and\/ $X\sq\dom P$ 
is\/ $\is11$ then there is a\/ $\id11$ set\/ $Y\sq\dom P$ 
and a\/ $\id11$ function\/ $F : Y\to D$ such that\/ $X\sq Y$ 
and\/ $F\sq P$.
\ele
\bpf
The set $X_0=\dom P$ is $\ip11$ since $\ip11$ is closed under 
$\sus y\in\id11$. 
Therefore by Separation there is a $\id11$ set $Y\yt X\sq Y\sq X_0$.
By $\ip11$ Uniformization, there is a $\ip11$ set $F\sq P$ such 
that $\dom F=Y$ and $Y$ is a function. 
To show that $F$ is in fact $\id11$ note that 
$F(x)=y$ iff $x\in Y$ and 
$\kaz y'\in D\,(y\ne y'\imp\ang{x,n'}\nin F)$, 
which leads to a $\is11$ definition.
\epf

\vyk{

\subsection{Near-counterexamples}
\label{sec2}

The following examples show that  
Theorem \ref{mt}\ref{mt1} is not true any more 
for different extensions of the domain of 
$\fs11$ suborders of a Borel partial quasi-orders, 
such as $\fs11$ and $\fp11$ linear quasi-orders --- 
not necessarily suborders of Borel orderings, as well as 
$\fd12$ and $\fp11$ suborders of Borel orderings.
In each of these classes, a counterexample of cofinality 
$\omi$ will be defined.

\bex
[$\fs11$ LQO] 
\lam{ex1}
Consider a recursive coding of sets of rationals by reals. 
Let $Q_x$ be the set coded by a real $x$.
Let $X_\al$ be the set of all reals $x$ such that the 
maximal well-ordered initial segment
of $Q_x$ has the order type $\al$. 
We define 
$$
x\le y
\quad\text{iff}\quad
\sus\al\,\sus\ba\:(x\in X_\al\land y\in X_\ba\land\al\le\ba).
$$
Then $\le$ is a $\is11$ LQO on $\bn$ of cofinality $\om_1$.

Note that the associated strict order $x<y$,
iff $x\le y$ but not $y\le x$, 
is then more complicated than just $\fs11$,
therefore there is no contradiction in this example 
to the result mentioned in Remark~\ref{fn2}.\qed
\eex

\bex
[$\fp11$ LQO] 
\lam{ex2}
Let $D\sq\bn$ be the $\ip11$ set of codes of (countable) ordinals. 
Then
$$
x\le y
\quad\text{iff}\quad
x,y\in D\,\land\, |x|\le |y|
$$ 
is a $\ip11$ LQO of cofinality $\omi$.
Note that $\le$ is defined on a non-Borel $\ip11$ set $D$, 
and there is no $\fp11$ LQO of cofinality $\omi$ but 
defined on a Borel set --- by exactly the same argument as 
in Remark~\ref{fn2}.
\qed
\eex

\bex
[$\fp11$ LO] 
\lam{ex3}
To sharpen Example \ref{ex2}, define
$$
x\le y
\quad\text{iff}\quad
x,y\in D\;\land\;
\big( {|x|<|y|}\,\lor\, 
{(|x|=|y|\land x\lex y)}\big);
$$ 
this is a $\ip11$ LO of cofinality $\omi$.\qed
\eex

\bex
[$\fd12$ suborders] 
\lam{ex4}
Let $\le$ be the eventual domination order on $\om^\om$. 
Assuming the axiom of constructibility $\rV=\rL$, one 
can define a 
strictly $\le$-increasing $\id12$ $\omi$-sequence 
$\sis{x_\al}{\al<\omi}$ in $\om^\om$.\qed 
\eex

\bex
[$\fp11$ suborders] 
\lam{ex5}
Define a PQO $\le$ on $(\om\bez\ans0)^\om$ so that
$$
x\le y
\quad\text{ iff }\quad
\text{either}\quad x=y \quad\text{or}\quad
\lim_{n\to\iy}\:y(n)\,{/\hspace{-1ex}/}\,x(n)=\iy 
$$
(the ``or'' option defines the associated strict order $<$).
Assuming the axiom of constructibility $\rV=\rL$, define a 
strictly increasing $\id12$ $\omi$-sequence 
$\sis{x_\al}{\al<\omi}$ in $\om^\om$. 
By the Novikov -- Kondo -- Addison $\ip11$ Uniformization 
theorem, there is a $\ip11$ set 
$\sis{\ang{x_\al,y_\al}}{\al<\omi}\sq\bn\ti\dn$. 
%
%
%
Let $z_\al(n)=3^{x_\al(n)}\cdot 2^{y_\al(n)}$, $\kaz n$.
Then the \ddd\omi sequence $\sis{z_\al}{\al<\omi}$ is $\ip11$ 
and strictly increasing: indeed, factors of the form 
$2^{y_\al(n)}$ are equal 1 or 2 whenever $\al\in\dn$.\qed
\eex

}

\subsection{Ingredient 1: coding $\id11$ functions} 
\las{FF}

The proof of Theorem~\ref{mt} involves several technical methods 
of rather general nature, which we present in the three 
following sections.

Recall that $\ck$ is the 
least non-recursive (= the least non-$\id11$) ordinal.     
If $\al<\ck$ then let $\cf\al$ be the set of all $\id11$ 
\lr\ order preserving maps $F:\stk{\bn}\cle\to\stk{2^\al}\lexe$, 
so that
\imar{F1}
\bus
\label{F1}
x\cle y\limp F(x)\lexe F(y)\quad\text{for all}\quad 
x,y\in\bn\,.
\eus
Such a function $F$ has to be \dd\apr invariant on $\bn$.
Let $\cF=\bigcup_{\al<\ck}\cf\al$.

If, in addition, $X\sq\bn$ is a $\is11$ set then let $\cfx X$ 
consist of 
all $\id11$ functions $F\in\cF$ such that  
\imar{F2}
\bus
\label{F2}
x,y\in X\,\text{ are \dd\cle incomparable} \limp  F(x)=F(y)
\,.
\eus
\vyk{
or equivalently, 
$$
\kaz x,y\in X\:
\left\{
\bay{rll}
x\cle y&\imp& F(x)\lexe F(y)\,,\quad\text{and}\\[1ex]
x,y\,\text{ are \dd\cle incomparable}&\imp& F(x)=F(y)\,.
\eay
\right.
$$
}%
or equivalently, $F(x)\lex F(y)\imp x\cl y$ for all $x,y\in X$. 

Note that a function $F\in\cfx X$ has to be not 
just \dd\apr invariant, 
but also invariant \vrt\ the common equivalence hull of the 
relation $\apr$ 
and the (non-equivalence) relation of being \dd\cle incomparable. 
In particular, if for any $x,y\in X$ 
there is $z\in X$ \dd\cle incomparable with both $x$ and $y$, 
then the only maps in $\cfx X$ 
are those constant on $X$. 

\bdf
\lam{dEF}
Let, for  $x,y\in\bn$:
$\,x\rEF y\;$ iff $\;\kaz F\in\cF\:(F(x)=F(y))$,
$$
\hspace*{-10ex}
x\refx X y \;\;\text{ iff }\;\; \kaz F\in\cfx X\:(F(x)=F(y))\,\;\;
\text{(here $X\sq\bn$ is $\is11$)}.
\hspace*{-6ex}
\eqno\qed 
$$
\edf

\ble
\lam{Lef}
$\rEF$ is a smooth\/ $\is11$ \eqr, and 
if\/ $R(x,y)$ is a\/ $\ip11$ relation and\/ 
$\kaz x,y\:({x\rEF y}\imp R(x,y))$ 
then there is a single function\/ $F\in\cF$ such that\/ 
$\kaz x,y\:({F(x)=F(y)}\imp R(x,y))$. 

Similarly, if\/ $X\sq\bn$ is a\/ $\is11$ set then\/ 
$\refx X$ is a smooth\/ $\is11$ \eqr, and 
if\/ $R(x,y)$ is a\/ $\ip11$ relation and\/ 
$\kaz x,y\:({x\refx X y}\imp R(x,y))$ 
then there is a function\/ $F\in\cfx X$ such that\/ 
$\kaz x,y\:({F(x)=F(y)}\imp R(x,y))$. 
\ele
\bpf
We concentrate on the first part; the result for the second part 
is pretty similar. 
We'll make use of an appropriate 
coding of functions in $\cF$, based on  
a standard coding system of $\id11$ sets.
A \rit{code} will be a such-and-such pair 
$f=\ang{\ve,k}\in\bn\ti\om$.  
We require that:
\ben
\Renu
\itla{cod1}
the relation $\le_\ve=\ens{\ang{i,j}}{\ve(2^i\cdot3^j)=0}$ is a 
(non-strict) wellordering of the set 
$\dom{(\le_\ve)}$ ---
in this case, we let: 
\bit
\item[$-$]\msur 
$|\ve|=\otp\ve<\omi$ be the order type of $\le_\ve$,
 
\item[$-$]\msur 
$\ba_\ve:\dom{(\le_\ve)}\onto|\ve|$ be the order-preserving 
bijection, 

\item[$-$]\msur 
$H_\ve:\bn\onto(\bn){}^{|\ve|}$ be the induced homeomorphism;    
\eit

\itla{cod2}
$k$ belongs to the set $\dB\sq\om$ of 
\rit{codes of\/ $\id11$ sets} 
$B\sq\bn\ti\bn$, 
so that it is assumed that $\dB$ is a $\ip11$ set, and for any 
$k\in\dB$ a $\id11$ set $B_k\sq\bn\ti\bn$ is defined, and 
conversely, for any $\id11$ set $B\sq\bn\ti\bn$ there is a code 
$k\in\dB$ with $B=B_k$, and finally there exist two $\ip11$ 
sets $W,W'\sq\om\ti\bn\ti\bn$ such that if $k\in\dB$ and 
$x,y\in\bn$ then
$$
\ang{x,y}\in B_k\leqv W(k,x,y)\leqv\neg\:W'(k,x,y)\,; 
$$

\itla{cod3}
we define 
$\fun= \ens{k\in\dB}{B_k\text{ is a total map }\bn\to\bn}$, 
the set of codes of all $\id11$ functions $F:\bn\to\bn$ --- 
this is still a $\ip11$ set because the key condition 
$\dom B_f=\bn$ can be expressed by 
$$
\kaz x\:\sus y\in\id11(x)\:W(k,x,y)\,,
$$
where the quantifier $\sus y\in\id11(x)$ is known to preserve 
the type $\ip11$;

\itla{cod4}
if $\ve$ satisfies \ref{cod1} and $k\in\fun$ then let 
$F^\ve_k$ be the $\id11$ map $\bn\to(\bn){}^{|\ve|}$ defined by 
$F^\ve_k(x)=H_\ve(B_k(x))$ 
for all $x\in\bn$;
\een

\bdf
\lam{ff}
Let $\dF$ be the set of all pairs $f=\ang{\ve,k}$ such that 
$\ve\in\bn$ satisfies \ref{cod1}, $\ve\in\id11$, $k\in\fun$, and 
$F^\ve_k\in\cf{|\ve|}$.

If $X\sq\bn$ is a $\is11$ set then let 
$\dF(X)=\ens{\ang{\ve,k}\in\dF}{F^\ve_k\in\cfd{|\ve|}X}$.\qed
\edf

The following is a routine fact.

\bcl
\lam{31.1}
$\dF\sq\bn\ti\om$ is a countable\/ $\ip11$ set of\/ 
$\id11$ elements, and\/ 
$\cF=\ens{F^\ve_k}{\ang{\ve,k}\in\dF}$. 
If $X\sq\bn$ is a $\is11$ set then $\dF(X)\sq\bn\ti\om$ is a 
countable\/ $\ip11$ set of\/ 
$\id11$ elements, and\/ 
$\cfx X=\ens{F^\ve_k}{\ang{\ve,k}\in\dF(X)}$. \qed
\ecl

In continuation of the proof of Lemma~\ref{Lef}, note that 
$$
x\rEF y
\eqv
\kaz \ang{\ve,k}\in\dF\:(F^\ve_k(x)=F^\ve_k(y)) 
\eqv
\kaz \ang{\ve,k}\in\dF\:(B_k(x)=B_k(y)), 
$$
and this easily implies that $\rEF$ is $\is11$ by Claim~\ref{31.1}.
Now prove the claim of Lemma~\ref{Lef} related to  $R(x,y)$.  
We re-write the assumption as follows:
$$
\kaz x,y\:(\neg\:R(x,y)\limp \neg\:(x\rEF y))\,,
$$
or, equivalently by Claim~\ref{31.1}, as 
$$
\kaz x,y\:\big(\neg\:R(x,y)\limp 
\sus\ang{\ve,k}\in\id11\:
\underbrace
{(\ang{\ve,k}\in\dF\land F^\ve_k(x)\ne F^\ve_k(y))}
_{P(x,y\,;\,\ve,k)}\;
\big)\,.
$$
The relation $P$ is expressible by a $\ip11$ formula 
by means of \ref{cod2} and Claim~\ref{31.1}.
It follows by 
Lemma~\ref{ks} 
that 
there is a $\id11$ set $W\sq\bn\ti\bn$ satisfying
$\neg\:R(x,y)\limp W(x,y)$,
and a $\id11$ map $\Phi(x,y)=\ang{\ve(x,y),k(x,y)}:W\to\dF$ 
such that we have 
$F^{\ve(x,y)}_{k(x,y)}(x)\ne F^{\ve(x,y)}_{k(x,y)}(y)$  
for all $\ang{x,y}\in W$ 
--- 
then, in particular, for all $x,y$  
with $\neg\:R(x,y)$. 

The range $H=\ens{\Phi(x,y)}{\ang{x,y}\in W}$ is then a 
$\is11$ subset of the (countable) $\ip11$ set $\dF$. 
By Separation, there is a $\id11$ set $D$ with $H\sq D\sq\dF$.
As a countable $\id11$ set, it admits a $\id11$ enumeration 
$D=\ens{\ang{\ve_n,k_n}}{n\in\dN}$, and by construction 
$\kaz n\,(F^{\ve_n}_{k_n}(x)=F^{\ve_n}_{k_n}(y))$ implies 
$R(x,y)$.
Let 
$$
F(x)=
F^{\ve_0}_{k_0}(x)\we F^{\ve_1}_{k_1}(x)\we F^{\ve_2}_{k_2}(x)\we\dots
$$ 
for $x\in\bn.$ 
Then $F\in\cF$ and 
$F(x)=F(y)\imp R(x,y)$.
\epF{Lemma~\ref{Lef}}

\subsection{Ingredient 2: invariant separation}
\las{Ti}

In the assumptions of Theorem~\ref{mt}, let $\rE$ be a $\is11$ 
\eqr\ containing $\apr$ 
(so that $x\apr y$ implies $x\rE y$). 
A set $X\sq\bn$ is 
\rit{downwards\/ \dd\cle closed in each\/ \dd\rE class} 
iff we have ${x\in X}\imp {y\in X}$ whenever $x\rE y$ and $y\cle x$. 
The notion of a set 
\rit{upwards\/ \dd\cle closed in each\/ \dd\rE class} 
is similar. 

\ble
\lam{i1}
Let\/ $\rE$ be a $\is11$ \eqr\ containing $\apr$, 
$X,Y$ be disjoint\/ $\is11$ sets, 
satisfying\/ $y\ncle x$ whenever\/ $x\in X\land y\in Y\land x\rE y$. 
Then there is a\/ $\id11$ set\/ $Z$, 
downwards\/ \dd\cle closed in each\/ \dd\rE class and satisfying\/ 
$X\sq Z$ and\/ $Y\cap Z=\pu$.
\vyk{
If moreover\/ $x\cl y$ holds  
whenever\/ $x\in X\land y\in Y\land x\rE y$, 
then we can additionally require that\/ $x\cl y$ holds  
whenever\/ $x\in Z\land y\in Y\land x\rE y$.
}%
\ele
\bpf
Let $Y'=\ens{y'}{\sus y\in Y\,(y\cle y')}$; 
still $Y'\cap X=\pu$ and $Y'$ is $\is11$.
Using Separation, define an increasing sequence of sets
$$
X = X_0\sq A_0\sq X_1\sq A_1\sq\dots\sq X_n\sq A_n\sq\dots
\sq \bn\bez Y'
$$
so that $A_{n}\in \id11$ and  
$X_{n+1} = \ens{x'\in\bn}
{\sus x \in A_n\,({x'\rE x}\land x'\cle x)}$ 
for all $n$. 
If $A_n \cap Y'=\pu$ then $X_{n+1}\cap Y'=\pu$ as well since 
$Y'$ is upwards closed, which 
justifies the inductive construction. 
Furthermore, a proper execution of the construction 
yields the final set 
$Z= \bigcup_n A_n =  \bigcup_nX_n$ in $\id11$.
(We refer to the proof of an ``invariant'' effective separation 
theorem in \cite{hkl} or a similar construction 
in \cite[Lemma 10.4.2]{kanB}.) 
Note that by construction $X\sq Z$, but $Z\cap Y =\pu$, and 
$Z$ is downwards \dd\cle closed in each \dd{\rE}class.
%
\epf

\bcor
\lam{i2}
Let\/ $\rE$ be\/ $\rE_\cF$. 
If\/ $X,Y\sq\bn$ are disjoint\/ $\is11$ sets and\/ 
$\eke X\cap\eke Y\ne\pu$ then there are points\/
$x\in X$, $y\in Y$ with\/ $x\rE y$ and\/ $y\cle x$.
\ecor
\bpf 
Otherwise by Lemma~\ref{i1} there is a $\id11$ set $Z$ such that 
$X\sq Z$ and $Y\cap Z=\pu$, 
and  downwards\/ \dd\cle closed in each\/ \dd\rE class.  
Then, by Lemma~\ref{Lef}, there is a function $F\in\cF$ such that 
$x\in Z\imp y\in Z$ holds whenever $F(x)=F(y)$ and $x\cle y$.
It follows that  the derived  function
$$
G(x)=\left\{
\bay{rcl}
F(x)\we 0\,,&\text{whenewer}& x\in Z\\[1ex]
F(x)\we 1\,,&\text{whenewer}& x\in\bn\bez Z 
\eay
\right.
$$
belongs to $\cF.$ 
Thus if $x\in Z$ and $y\nin Z$, say $x\in X$ and $y\in Y$, 
then $G(x)\ne G(y)$ and hence $x\nE y$, a contradiction. 
\epf

\vyk{
\bcor
\lam{i3}
Let\/ $\pu\ne A\sq\bn$ is a\/ $\is11$ set and\/ $\rE$ be\/ 
$\refx A$. 
If\/ $X,Y\sq A$ are disjoint\/ $\is11$ sets and\/ 
$\eke X\cap\eke Y\ne\pu$ then there are points\/
$x\in X$ and\/ $y\in Y$ satisfying\/ $x\rE y$ and\/ $x\ncl y$.
\ecor
\bpf
Similar to the proof of Corollary~\ref{i2}, yet one has to refer to 
the the moreover claim of Lemma~\ref{i1} in order to check that 
$G$ belongs to $\cfx A$.
\epf
}

\subsection{Ingredient 3: the Gandy -- Harrington forcing} 
\las{0gh}

The Gandy -- Harrington forcing notion $\gh$ 
is the set of all $\is11$ sets $\pu\ne X\sq \bn,$
ordered so that smaller sets are stronger conditions.
We also define $\gh_n$ ($n\ge2$) to be the set of all 
$\is11$ sets $\pu\ne X\sq (\bn){}^n$.

It is known that $\gh$ adds a point of $\bn,$ whose 
name will be $\dox$.

\vyk{
\bre
\lam{notr}
If $x\in\bn$ is $\id11$ then the set $X=\ans x$ belongs to 
$\gh$, of course. 
To get rid of such trivial conditions, they consider 
\rit{the domain of nontriviality} of all conditions $X\in\gh$ 
not containing $\id11$ elements, and hence uncountable.
\ere
}

Together with $\gh$, some other related forcing notions 
will be considered below, for instance, the product 
$\ghd=\gh\ti\gh$ which consists of all cartesian 
products of the form $X\ti Y,$ where $X,Y\in\gh.$ 
It follows from the above
that $\ghd$ forces a pair of points of $\dn,$ whose 
name will be $\ang{\doxl,\doxr}$.

There is another important subforcing introduced in \cite{hms}.
If $\rE$ is a $\is11$ \eqr\ on $\bn$ then 
let $\dpe$ consist of all sets of the form 
$X\ti Y,$ where $X\yi Y\in\gh$ and 
$(X\ti Y)\cap{\rE}\ne\pu$.

A condition $X\ti Y$ in $\dpe$ is \rit{saturated} iff 
$\eke X=\eke Y$.

\ble
\lam{r6}
If\/ $X\ti Y$ is a condition in $\dpe$ then there is a 
stronger saturated subcondition\/ $X'\ti Y'$ in $\dpe$.
\ele
\bpf  
$X'=\ens{x\in X}{\sus y\in Y\,(x\rE y)}$, 
$Y'=\ens{y\in Y}{\sus x\in X\,(x\rE y)}$.
\epf

\bre
\lam{r7}
If $X\ti Y$ is a 
saturated condition in $\dpe$, and 
$\pu\ne X'\sq X$ is a $\is11$ set, then $Y'=Y\cap\eke{X'}$ is 
$\is11$ and $X'\ti Y'$ is still a saturated condition in $\dpe.$ 
It follows that $\dpe$ forces a pair of 
\dd\gh generic reals,  whose names 
will be $\doxl$ and $\doxr$ as above.
\ere


\ble
[2.9 in \cite{hms}]
\lam{Esmu}
Suppose that\/ $\rE$ is a smooth\/ $\is11$ \eqr.
Then\/ $\dpe$ forces\/ $\doxl\rE\doxr$.\qed 
\ele

Note that Lemma~\ref{Esmu}, generally speaking, fails in the 
non-smooth case. 

The next result will be pretty important.

\ble
[2.9 in \cite{hms}]
\lam{29}
Suppose that\/ 
$\cle$ is a\/ $\id11$ PQO on\/ $\bn,$
and for any\/ $A\in\gh$ there is a\/ $\is11$ \eqr\ $\rE_A$ 
on\/ $\bn$ such that if\/ $A\sq B$ then\/ 
$x\rE_A y$ implies\/ $x\rE_B y$.
Assume that\/ $X^\ast\in\gh$, and if\/ $B\in\gh\yt B\sq X^\ast$ 
then\/ $B\ti B$ does\/ 
{\ubf not} 
\dd{(\drof\gh{\rE_B})}force
that\/ $\doxl,\doxr$ are\/ \dd\cle comparable. 

Then\/ $\Xa$ is not\/ \dd\cle thin, in other words, 
there is a perfect set\/ $Y\sq X^\ast$ of pairwise\/ 
\dd\cle incomparable elements.\qed 
\ele

\vyk{
\bpf
Let $T$ be the set of all finite trees $t\sq2\lom.$
If $t\in T$ then let $M(t)$ be the set of all \dd\su maximal 
elements of $t$. 
Let $\Phi$ be the set of systems $\xi=\sis{X_u}{u\in t}$ of sets 
$X_u\in\gh\yt X_u\sq X^\ast$, such that $t\in T$ and 
\ben
\renu
\itla{291}
$X_\La = X^\ast$ (where $\La$ is the empty string);

\itla{292}
if $u\su v\in t$ then $X_v\sq X_u$;

\itla{293} 
if $u\we0$ and $u\we 1$ belong to $t$ then 
$X_{u\we 0}\ti X_{u\we1}$ belongs to $\drof{\gh}{\rE_{P_u}}$ 
and \dd{(\drof{\gh}{\rE_{P_u}})}forces that 
$\doxl$ is \dd\cle incomparable to $\doxr$;

\itla{294} 
there is a sequence $\sis{x_u}{u\in M(t)}$ of points $x_u\in X_u$ 
such that if $u,v\in M(t)$ then $x_u\rE_{X_{u\land v}}x_v$, 
where $u\land v$ is the largest string $w\in2\lom$ such that 
$w\su u$ and $w\su v$.
\een 
Say that a system $\sis{X_u}{u\in t}\in\Phi$ is 
\rit{saturated} if in addition
\ben
\renu
\atc\atc\atc\atc
\itla{295}
for any $v\in M(t)$ and $x\in X_v$ there is a sequence 
$\sis{x_u}{u\in M(t)}$ as in \ref{294}, such that $x_u=x$. 
\een

\bct
\lam{29d}
For any system\/ $\sis{X_u}{u\in t}\in\Phi$ there is 
a saturated system\/ $\sis{X'_u}{u\in t}$ in $\Phi$ such that 
$X'_u\sq X_u$ for all\/ $u$ and\/ $X'_u=X_u$ for all\/ 
$u\in t\bez M(t)$.
\ect
\bpf[claim]
If $u\in M(t)$ then simply let $X'_u$ be the set of all points 
$x\in X_u$ such that $x=x_u$ for some sequence 
$\sis{x_u}{u\in M(t)}$ as in \ref{294}, while if $u\in t\bez M(t)$ 
then keep $X'_u=X_u$.
\epf

\bct
\lam{29c}
For any\/ $t\in T$ and a saturated system\/ 
$\xi=\sis{X'_u}{u\in t}\in\Phi$, 
if\/ $u\in M(t)$ then there are sets $X_{u\we0}\yi X_{u\we1}$ 
such that the system $\xi$ extended 
by those sets still belongs to $\Phi$.
\ect
\bpf[claim]
As $Y_u\in\gh$ and $Y_u\sq X^\ast$, the condition 
$Y_u\ti Y_u$ does\/ 
{\ubf not} 
\dd{(\drof\gh{\rE_{Y_u}})}force
that\/ $\doxl,\doxr$ are\/ \dd\cle comparable. 
Pick a condition $U\ti V$ in $\drof\gh{\rE_{Y_u}}$ which 
\dd{(\drof\gh{\rE_{Y_u}})}forces the opposite, \ie,
that $\doxl,\doxr$ are \dd\cle incomparable.
By Lemma~\ref{le1}, we may assume that $U\ti V$ is 
\dd{\rE_{Y_u}}saturated.
We assert that the sets $X_{u\we0}=U$ and $X_{u\we1}=V$ prove 
the claim. 
Indeed it's enough to check \ref{294} for the extended system. 

Fix any $x\in X_{u\we0}=U$.
Then $x\in X_{u}$, therefore, as the given system is saturated,
there is a sequence $\sis{x_v}{v\in M(t)}$ of points 
$x_v\in X_v$ as in \ref{294}, such that $x_u=x$. 
On the other hand, as $U\ti V\in \drof\gh{\rE_{X_u}}$ 
is saturated, there is 
a point $y\in V=X_{u\we1}$ such that $x\rE_{X_u}y$.
We put $x_{u\we0}=x$ and $x_{u\we1}=y$.
\epf

We come back to the proof of Lemma~\ref{29}.
We let $t_n=2^{\le n}$ (all binary strings of length $\le n$)
Define, by induction, saturated systems 
$\xi_n=\sis{X_u}{u\in t_n}\in\Phi$
so that each $\xi_{n+1}$ extends $\xi_n$ by one layer of sets
$X_{u\we i}$, $u\in 2^{n}$ and $i=0,1$, obtained by
consecutive $2^n$ splitting operations as in Claim \ref{29c},
followed by the saturating reduction as in Claim \ref{29d}.

To introduce a genericity condition, 
fix a cardinal $\ka$ large enough for $\bV_\ka$ to be 
an elementary submodel of the universe \vrt\ all properties 
(including forcing properties) of $\gh$ and related forcings, which 
matter in the proof of the lemma. 
Let $\gM$ be a countable elementary submodel of $\bV_\ka$, 
and $\cM$ be the Mostowski collapse of $\gM$ 
(a countable transitive model of a fragment of $\ZFC$). 
We require that the construction of systems $\xi_n$ viewed as a path 
through the tree of all possible constructions of this type is 
generic over $\gM$. 

Then for any $x\in \bn$ the intersection $\bigcap_mX_{x\res m}$ 
is a singleton, say $\ans{f(x)}$, where $f(x)\in\bn$ is a real 
\dd\gh generic over $m$, and $f:\bn\to\bn$ is continuous. 
Moreover, if $x\ne y\in\bn$ and $n$ is the largest number such 
that $x\res n=y\res n$ and $x(n)=0$, $y(n)=1$, then 
the pair $\ang{f(x),f(y)}$ is 
\dd{\drof{\gh}{\rE_{P_u}}}generic over $\cM$ --- 
where $u=x\res n=y\res n$ --- and 
\dd\cle incomparable to $\doxr$ by \ref{293}; 
in particular, $f(x)\ne f(y)$. 
Therefore $Y =\ens{f(x)}{x\in\bn}\sq X$ is a perfect set 
of mutually 
incomparable reals.
\epf
}

The forcing $\gh$, as well as some of its derivates like $\dpe$, 
will be used below as forcing notions over the ground set 
universe $\bV$. 

\ble
[see \cite{hkl,hms}]
\lam{xx}
If\/ $X\in\gh$ then\/ $X$ \dd\gh forces that\/ $\dox\in X.$ 
Moreover if\/ $\Phi(x)$ is a\/ $\ip12$ formula and\/ $\Phi(x)$ 
holds for all\/ $x\in X$ then\/ $X$ \dd\gh forces that\/ 
$\dox$ satisfies\/ $\Phi(X)$. 

The same is true for other similar forcing notions like\/ 
$\dpe$.\qed 
\ele

Here (and below in some cases), given a $\is11$ (or $\ip11$) set 
$X$ in the ground universe $\bV$, 
{\ubf we denote by the same letter $X$ the extended set} 
(\ie, defined by the same formula) 
{\ubf in any generic extension of $\bV$}. 
By the Shoenfield absoluteness theorem, there is no ambiguity here. 
See \cite[2.4]{ksz} in more detail.

\subsection{Bounding thin partial orderings}
\las{s3}

Here {\ubf we prove claim \ref{mtb} of Theorem~\ref{mt}}. 
We'll make use of the family $\cF$ of 
\pagebreak[0]
$\id11$ functions, introduced in Section~\ref{FF}, and the 
corresponding smooth $\is11$ \eqr\ ${\rE}={\rEF}$.
Then $\approx$ is a subrelation of $\rE$ by Lemma~\ref{Lef}. 

The following partition on cases is quite common in this sort of 
proofs.\vom

{\ubf Case 1:} 
$\approx$ and $\rE$ coincide on  $\Xa$, so that 
${x\rE y}\imp{x\apr y}$ for $x,y\in\Xa$. 
Then, by Lemma~\ref{Lef}, there is a single function $F\in\cF$ 
such that $F(x)=F(y)$ implies $x\apr y$ for all $x,y\in\Xa$, 
as required. 
\vtm

{\ubf Case 2:} 
$\approx$ is a \rit{proper} subrelation of $\rE$ on $\Xa$, hence, 
the $\is11$ set 
$$
\Va=\ens{x\in\Xa}{\sus y\in \Xa\:(x\not\approx y\land x\rE y)}
$$ 
is non-empty. 
Our final goal will be to infer a contradiction;   
then the result for Case 1 proves Claim \ref{mtb} 
of the theorem.

Note that  
$\Va\ti\Va$ is a saturated condition in $\dpd\rE$.

\ble
\lam{31+}
Condition\/ $\Va\ti\Va$ 
\dd{(\dpd{\rE})}forces that\/ $\doxl$ and\/ $\doxr$ are\/ 
\dd\cle incomparable. 
\ele
\bpf
Suppose to the contrary that a subcondition $Y\ti Z$ 
either forces $\doxl\apr\doxr$ or forces $\doxl\cl\doxr$. 
We will get a contradiction in both cases. 
Note that $Y,Z\sq\Va$ are non-empty $\is11$ sets and 
$\eke Y\cap\eke Z\ne\pu$.\vom

{\ubf Case A:} $Y\ti Z$ forces $\doxl\apr\doxr$.\vom

{\ubf Subcase A1:} the $\is11$ set 
$W=\ens{\ang{y,y'}\in Y\ti Y}{y\rE y'\land y'\napr y}$ 
is empty, or in other words ${\rE}$ coincides with ${\apr}$ on $Y.$
By the non-emptiness of $\Va$ at least one of the $\is11$ sets 
$$
B =\ens{x}{\sus y\in Y\,(x\rE y\land x\ncle y)}\,,\; 
B'=\ens{x}{\sus y\in Y\,(x\rE y\land y\ncle x)}
$$ 
\vyk{
$$
\bay{rcl}
B &=& \ens{x\in\Va}{\sus y\in Y\,(x\rE y\land x\ncle y)}\,,\\[1ex] 
B'&=&\ens{x\in\Va}{\sus y\in Y\,(x\rE y\land y\ncle x)}\,,
\eay
$$ 
}%
is non-empty; assume that, say, $B\ne\pu$. 
Consider the $\is11$ set 
$$
A=\ens{x}{\sus y\in Y\,(x\rE y\land x\cle y)}\,;
\quad Y\sq A\,.
$$
Then $A\cap B=\pu$, 
$A$ is downwards closed while $B$ is 
upwards closed in each \dd\rE class, therefore $y\ncle x$ 
whenever $x\in A$, $y\in B$, and $x\rE y$.
Then $\eke A\cap\eke B=\pu$ by Corollary~\ref{i2}.
Yet by definition $\eke Y\cap\eke B\ne\pu$ and $Y\sq A$, 
which is a contradiction.\vom

{\ubf Subcase A2:} $W\ne\pu$.
Then the forcing notion $\gh(W)$ of all non-empty $\is11$ 
sets $P\sq W$ adds pairs of \dd\gh generic (separately) 
reals $y,y'\in Y$ which belong to $W$  
and satisfy $y'\rE y$ and $y'\napr y$, by Lemma~\ref{xx}. 

If $P\in\dpw$ then obviously $\eke{\dom P}=\eke{\ran P}$. 

Consider a more complex forcing notion
$\cP=\dpw\ti_{\rE}\gh$ of all pairs $P\ti Z'$, where $P\in\dpw$,  
$Z'\in\gh$, $Z'\sq Z$, and $\eke{\dom P}\cap\eke{Z'}\ne\pu$. 
For instance, $W\ti Z\in\dpw\ti_{\rE}\gh$.  
Then  
$\cP$ adds a pair $\ang{\doxl,\doxr}\in W$ and a separate 
real $\dox\in B$ such that both pairs $\ang{\doxl,\dox}$ and 
$\ang{\doxr,\dox}$ are \dd{(\dpx{})}generic, hence, we have 
$\doxl\apr\dox\apr\doxl$ 
(in the extended universe $\bV[\doxl,\doxr,\dox]$) 
by the choice of $Y\ti Z$. 
On the other hand, $\doxl\napr\doxl$ by Lemma~\ref{xx}, 
since the pair belongs 
to $W$, which is a contradiction.\vom

{\ubf Case B:} $Y\ti Z$ forces $\doxl\cl\doxr$.\vom

{\ubf Subcase B1:} the $\is11$ set 
$ 
W=\ens{\ang{y,z}\in Y\ti Z}{{z\rE y}\land {z\cle y}}
$ 
is empty. 
Then the $\is11$ sets  
$$
Y_0=\ens{y'}{\sus y\in Y\,(y\rE y'\land y'\cle y)}\,,\;\;
Z_0=\ens{z'}{\sus z\in Z\,(z\rE z'\land z\cle z')} 
$$
are disjoint and \dd\cle closed resp.\ downwards and upwards, 
hence we have $\eke{Z_0}\cap\eke{Y_0}=\pu$ by Corollary~\ref{i2}. 
However $\eke{Z}\cap\eke{Y}\ne\pu$, 
which is a contradiction as $Z\sq Z_0$, $Y\sq Y_0$.\vom

{\ubf Subcase B2:} $W\ne\pu$.
Consider the forcing $\dpw$ of all non-empty $\is11$ sets 
$P\sq W$; if $P\in\dpw$ then obviously $\eke{\dom P}=\eke{\ran P}$. 
Consider a more complicated forcing $\dpwe$ of all products 
$P\ti Q$, where $P,Q\in\dpw$ and 
$
\eke{\dom P}\cap\eke{\dom Q}\ne\pu\,.
$
In particular $W\ti W\in \dpwe$.

Let $\ang{x,y;x',y'}$ be a \dd\dpwe generic quadruple in $W\ti W$, 
so that both $\ang{x,y}\in W$ and $\ang{x',y'}\in W$ are 
\dd\dpw generic pairs in $W$, 
and both $y\cle x$ and  $y'\cle x'$ hold by the definition 
of $W$.
On the other hand, an easy argument shows that both criss-cross 
pairs $\ang{x,y'}\in X\ti Y$ and $\ang{x',y}\in X\ti Y$ are 
\dd{\dpd{\rE}}generic, hence $x\cl y'$ and $x'\cl y$ by the 
choice of $X\ti Y$. 
Altogether  $y\cle x\cl y'\cle x'\cl y$, which is a 
contradiction.
\epf

\vyk
{

\ble
\lam{31-4}
If\/ $Y\sq\Va$ is a non-empty\/ $\is11$ set then 
there exist\/ $x,y\in Y$ such that\/ $x\rE y$ and\/ 
$x\napr y$.
\ele
\bpf 
Assume to the contrary that ${\rE}$ coincides with ${\apr}$ on $Y.$
Still by the non-emptiness of $\Va$ at least one of the $\is11$ sets 
$$
B = \ens{x}{\sus y\in Y\,(x\rE y\land x\ncle y)}\,,\;
B'=\ens{x}{\sus y\in Y\,(x\rE y\land y\ncle x)} 
$$ 
\vyk{
$$
\bay{rcl}
B &=& \ens{x\in\Xa}{\sus y\in Y\,(x\rE y\land x\ncle y)}\,,\\[1ex]
B'&=&\ens{x\in\Xa}{\sus y\in Y\,(x\rE y\land y\ncle x)}\,,
\eay
$$ 
}%
is non-empty; assume that, say, $B\ne\pu$. 
Consider the $\is11$ set 
$$
A=\ens{x}{\sus y\in Y\,(x\rE y\land x\cle y)}\,;
\quad Y\sq A\,.
$$
Then $A\cap B=\pu$ (since ${\rE}={\apr}$ on $Y$), 
$A$ is downwards \dd\cle closed while $B$ is 
upwards \dd\cle closed in each \dd\rE class --- therefore we have 
$y\ncle x$ whenever $x\in A$, $y\in B$, and $x\rE y$.
By Corollary~\ref{i2}, $\eke A\cap\eke B=\pu$.
Yet by definition $\eke Y\cap\eke B\ne\pu$ and $Y\sq A$, 
which is a contradiction. 
\epF{Lemma}

\ble
\lam{31e}
Condition\/ $\Va\ti\Va$ 
\dd{(\dpd{\rE})}forces that\/ $\doxl\napr\doxr$.
\ele
\bpf
Suppose to the contrary that a subcondition 
$Y\ti Z$ \dd{(\dpd\rE)}forces $\doxl\apr\doxr$. 
Then any \dd{(\dpd\rE)}generic pair $\ang{y,z}\in Y\ti Z$ 
satisfies $y\apr z$. 
It follows that any two \dd\gh generic reals $y,y'\in Y$ with 
$y\rE y'$ satisfy $y\apr y'$, too, because one can define 
then a single $z\in Z$ such that both pairs $\ang{y,z}$ and 
$\ang{y',z}$ are \dd{(\dpx{})}generic.

Now we claim that the $\is11$ set 
$W=\ens{\ang{y,y'}\in Y\ti Y}{y\rE y'\land y'\napr y}$ 
is empty. 
Indeed, otherwise forcing $\gh(W)$ of all non-empty $\is11$ 
sets $P\sq W$ adds pairs of \dd\gh generic (separately) 
reals $y,y'\in Y$ which belong to $W$, hence, satisfy 
$y'\rE y$ and $y'\napr y$, which contradicts the assumption.

Thus if $y,y'\in X$ then $y\rE y'\imp y'\apr y$, 
contrary to Lemma~\ref{31-4}.
\epf

\ble
\lam{31f}
$\Va\ti\Va$ 
\dd{(\dpd{\rE})}forces\/ $\doxl\ncl\doxr$ and\/ $\doxr\ncl\doxl$.
\ele
\bpf
Otherwise there is a condition $X\ti Y$ in $\dpd{\rE}$ with 
$X\cup Y\sq\Va,$ which forces $\doxl\cl\doxr$;   
$X,Y$ are non-empty $\is11$ sets satisfying 
$\ek{X}{\rE}\cap\ek{Y}{\rE}\ne\pu$. 
If the set 
$ 
W=\ens{\ang{x,y}\in X\ti Y}{{x\rE y}\land {y\cle x}}
$ 
is empty then
\pagebreak[0]  
$$
X_0=\ens{x'}{\sus x\in X\,(x\rE x'\land x'\cle x)}\,,\;\;
Y_0=\ens{y'}{\sus y\in Y\,(y\rE y'\land y\cle y')} 
$$
are disjoint and \dd\cle closed resp.\ downwards and upwards 
$\is11$ sets, 
hence $\eke{X_0}\cap\eke{Y_0}=\pu$ by Corollary~\ref{i2}. 
However $\eke{X}\cap\eke{Y}\ne\pu$, 
which is a contradiction as $X\sq X_0$, $Y\sq Y_0$.
It follows that $W\ne\pu$.
By Lemma~\ref{r6} we can assume that  $X\ti Y$ is saturated, 
so that $\ek{X}{\rE}=\ek{Y}{\rE}$. 

Consider the forcing $\dpw$ of all non-empty $\is11$ sets 
$P\sq W$; if $P\in\dpw$ then obviously 
$\eke{\dom P}=\eke{\ran P}$. 
Consider the reduced product forcing $\dpwe$ of all products 
$P\ti Q$, where $P,Q\in\dpw$ and 
$$
\eke{\dom P}=\eke{\ran P}=\eke{\dom Q}=\eke{\ran Q}\,.
$$
In particular $W\ti W\in \dpwe$.
Let $\ang{x,y;x',y'}$ be a \dd\dpwe generic quadruple in $W\ti W$, 
so that both $\ang{x,y}\in W$ and $\ang{x',y'}\in W$ are 
\dd\dpw generic pairs in $W$, 
and we have both $y\cle x$ and  $y'\cle x'$ by the definition 
of $W$.
On the other hand, an easy argument shows that both criss-cross 
pairs $\ang{x,y'}\in X\ti Y$ and $\ang{x',y}\in X\ti Y$ are 
\dd{\dpd{\rE}}generic, therefore $x\cl y'$ and $x'\cl y$ by the 
choice of $X\ti Y$. 

Thus $y\cle x\cl y'\cle x'\cl y$, which is a 
contradiction.
\epF{Lemma}

\bcor
\lam{l9} 
$\Va\ti\Va$ \dd{(\dpd{\rE})}forces\/ 
$\doxl\ncle\doxr$ and\/ $\doxr\ncle\doxl$.\qed
\ecor
}

To accomplish the proof of \ref{mtb} of Theorem~\ref{mt}, 
note that by Lemma~\ref{31+} and Lemma~\ref{29}  
(with ${\rE_A}={\rE}$ for all $A$)
there is a perfect 2wise 
\dd\approx inequivalent set, so $\cle$ is not thin, 
contrary to our assumptions.

\subsection{Decomposing thin partial orderings}
\las{Td}

We prove {\ubf claim \ref{mtd} of Theorem~\ref{mt}} in this Section. 
Let $\Ua$ be the $\is11$ set of all reals $x\in\Xa$ such that 
there is no $\id11$ \dd\cle chain $C$ containing $x$. 

{\ubf We assume to the contrary that $\Ua\ne\pu$}.

The proof will make heavy use of the functions in 
families of the form $\cfx X$, introduced in Section~\ref{FF}.
If $X\sq\bn$ is a $\is11$ set then ${\rE_X}={\refx X}$ 
is a smooth $\is11$ \eqr\ by Lemma~\ref{Lef}.

If $X\sq X'$ then $\cfx{X'}\sq\cfx X$, and hence 
$x\rE_X y$ implies $x\rE_{X'} y$.

\bcor
[of Lemma~\ref{Esmu}]
\lam{Esmu4}
If\/ $X\sq \Ua$ is a non-empty\/ $\is11$ set then the 
condition\/ $X\ti X$ 
\dd{(\dpx X)}forces that\/ $\doxl\rE_X\doxr$.\qed
\ecor

\ble
\lam{51+}
Let\/ $X\sq \Ua$ be a non-empty\/ $\is11$ set. 
Then\/ $X\ti X$ 
does {\bfit not}\/ \dd{(\dpx X)}force that\/ $\doxl,\doxr$ 
are\/ \dd\cle comparable.
\ele
\bpf
Suppose to the contrary that $X\ti X$ forces the comparability. 
Then there is a subcondition $Y\ti Z$ which
either forces $\doxl\apr\doxr$ or forces ${\doxl\cl\doxr}$; 
$Y,Z\sq X$ are non-empty $\is11$ sets and 
$\ek{Y}{\rE_X}\cap\ek{Z}{\rE_X}\ne\pu$.\vom

{\ubf Case A:} $Y\ti Z$ forces $\doxl\apr\doxr$.\vom

{\ubf Subcase A1:} the $\is11$ set 
$W=\ens{\ang{y,y'}\in Y\ti Y}{y\rE_X y'\land y'\napr y}$ 
is empty. 
Then $Y$ is a \dd\cle chain: 
indeed if $x,y\in Y$ are \dd\cle incomparable 
then by definition we have $x\rE_X y$, hence $x\apr y$, 
contradiction. 
Let $C$ be the $\ip11$ set of all reals \dd\cle comparable with 
each $y\in Y$; then $Y\sq C$. 
By Separation there is a $\id11$ set $D\yt Y\sq D\sq C$. 
Let $C'$ be the $\ip11$ set of all reals in $D$, 
\dd\cle comparable with each $d\in D$; then
$Y\sq C'\sq D$. 
Take any $\id11$ set $B$ with $Y\sq B\sq C'$. 
By construction $B$ is a $\id11$ \dd\cle chain with $\pu\ne Y\sq B$, 
contrary to the definition of $\Ua$. 
\vom

{\ubf Subcase A2:} $W\ne\pu$: yields a contradiction similarly 
to Subcase A2 in the proof of Lemma~\ref{31+}.\vom 

{\ubf Case B:} $Y\ti Z$ forces $\doxl\cl\doxr$.\vom

{\ubf Subcase B1:} the $\is11$ set 
$W=\ens{\ang{y,z}\in Y\ti Z}{y\rE_X z\land y\not\cl z}$ 
is non-empty. 
Let $Y'=\dom W$.
As $Y'\sq X,$ the condition $Y'\ti Y'$ \dd{(\dpx X)}forces that 
$\doxl,\doxr$ are \dd\cle comparable. 
Therefore by the result in Case A there is a condition 
$A\ti B$ in $\dpx X$, 
with $A\cup B\sq Y'$, which forces $\doxl\cl\doxr$; 
for if it forces $\doxr\cl\doxl$ then just consider $B\ti A$ 
instead of $A\ti B$.
Consider the forcing notion $\cP$ of all non-$\pu$ $\is11$ sets
of the form $P\ti B'$, where 
$$
P\sq W,\;\;\dom P\sq A,\;\; B'\sq B, 
\qand 
\ek{B'}{\rE_X}=\ek{\dom P}{\rE_X}=\ek{\ran P}{\rE_X}\,.
$$
For instance if $B'=B$ and $P=\ens{\ang{x,y}\in W}{x\in A}$ then 
$P\ti B'\in\cP$. 
Then $\cP$ forces a pair $\ang{\doxl,\doxr}\in W$ and a separate 
real $\dox\in B$ such that both pairs $\ang{\doxl,\dox}$ and 
 $\ang{\doxr,\dox}$ are \dd{(\dpx X)}generic. 
It follows that $\cP$ forces 
both $\doxl\cl\dox$ (as this pair belongs to $A\ti B$) 
and $\dox\cl\doxr$ (as this pair belongs to $Y\ti Z$), 
hence, forces  $\doxl\cl\doxr$. 
On the other hand $\cP$ forces $\doxl\not\cl\doxr$ 
(as this pair belongs to $W$), which is a contradiction.\vom

{\ubf Subcase B2:}  
$W=\pu$, in other words, if\/ $y\in Y,$ $z\in Z$, and\/ $y\rE_X z$ 
then\/ $y\cl z$ strictly.
Then by Lemma~\ref{i1}  
there is a\/ $\id11$ set\/ $\wY\sq\bn$, 
downwards\/ \dd\cle closed 
in each\/ \dd{\rE_X}class, such that $Y\sq\wY$ and  
still $Z\cap\wY=\pu$. 
We claim that, moreover, 
\bce
{\it if\/ $y\in \wY\cap X$, $z\in X\bez\wY$, and\/ $y\rE_X z$, 
then\/ $y\cl z$}. 
\ece
Indeed otherwise, the following $\is11$ set 
$$
H_0=\ens{z\in X\bez\wY}
{\sus y\in\wY\cap X\:(y\rE_X z\land y\not\cl z)}
\sq X
$$          
is non-$\pu$. 
As above, 
there is a saturated condition $H\ti H'$ in $\dpx X$, 
with $H\cup H'\sq H_0$, 
which forces $\doxl\cl \doxr$, and then  
$z\cl z'$ holds whenever $\ang{z,z'}\in H\ti H'$  and $z\rE_X z'$. 
By construction the $\is11$ set 
$$
\wY_1 =\ens{y\in\wY\cap X}
{\sus z'\in H'\,({y\rE_X z'} \land y\not\cl z')} 
$$
satisfies $\ek{\wY_1}{\rE_X}=\ek{H}{\rE_X}=\ek{H'}{\rE_X}$, 
hence $\wY_1\ti H$ is a condition in $\dpx X$.
Let $\ang{y_1,z}\in \wY_1\ti H$ be any \dd{(\dpx X)}generic pair. 
Then ${y_1\rE_X z}$ by Corollary~\ref{Esmu4}, 
and, by the choice of $X$ and the result in Case A, 
we have $y_1\cl z$ or $z\cl y_1$. 
However by construction $y_1\in\wY$, $z\nin \wY$, and $\wY$ 
is downwards closed in each \dd{\rE_X}class. 
Thus in fact $y_1\cl z$.
Therefore, for all $z'\in H'$, if $y_1\rE_X z'$ then 
$y_1\cl z\cl z'$, which contradicts to $y_1\in\wY_1$. 

Thus indeed $y\cl z$ holds whenever 
$y\in \wY\cap X$, $z\in X\bez\wY$, and\/ $y\rE_X z$.  
By Lemma~\ref{Lef}   
{\it there is a single function\/ $F\in\cF_X$ 
such that if\/ $y\in \wY\cap X$, $z\in X\bez\wY$, and\/ 
$F(y)=F(z)$, then\/ $y\cl z$}.

We claim that {\it the derived  function}
$$
G(x)=\left\{
\bay{rcl}
F(x)\we 0\,,&\text{whenewer}& x\in\wY\\[1ex]
F(x)\we 1\,,&\text{whenewer}& x\in\bn\bez\wY 
\eay
\right.
$$
{\it belongs to\/ $\cfx X.$} 
First of all, still $G\in\cF$ since $C$ is downwards 
\dd\cle closed in each \dd{\dpx X}class. 
Now suppose that $z,y\in X$ and $G(y)\lex G(z)$.
Then either $F(y)\lex F(z)$, or $F(z)=F(y)$ and $y\in\wY$ but 
$z\nin\wY$.
In the ``either'' case immediately $y\cl z$ since $F\in\cfx X$ 
\snos
{The family $\cF$ would not work in the passage; here  
we have to use $\cfx X$ instead.}
. 
In the ``or'' case we have $y\cl z$ by the choice of $F$ 
and the definition of $G$.
Thus $G\in\cfx X$.\vom 

Now pick any pair of reals $y\in Y$ and $z\in Z$ with $y\rE_X z$. 
Then we have $G(x)=G(y)$ since $G\in\cfx X$. 
But $y\in \wY$ and $z\nin\wY$ hold since 
$Y\sq\wY$ and $Z\cap\wY=\pu$ by construction, and in this case 
surely $G(y)\ne G(z)$ by the definition of $G$.
This contradiction completes the proof of Lemma~\ref{51+}.
\epf

\vyk{
\ble
\lam{51-4}
If\/ $Y\sq X\sq\Ua$ are non-empty\/ $\is11$ sets then 
there exist\/ $x,y\in Y$ such that\/ $x\rE_X y$ and\/ 
$x\napr y$.
\ele
\bpf 
Otherwise $Y$ is a \dd\cle chain. 
(Indeed if $x,y\in Y$ are \dd\cle incomparable 
then by definition we have $x\rE_X y$, hence $x\apr y$, 
contradiction.)  
Let $C$ be the $\ip11$ set of all reals \dd\cle comparable with 
each $y\in Y$; then $Y\sq C$. 
By Separation there is a $\id11$ set $D\yt Y\sq D\sq C$. 
Let $C'$ be the $\ip11$ set of all reals in $D$, 
\dd\cle comparable with each $d\in D$; then
$Y\sq C'\sq D$. 
Take any $\id11$ set $B$ with $Y\sq B\sq C'$. 
By construction $B$ is still a $\id11$ chain with $\pu\ne Y\sq B$, 
contrary to the definition of $\Ua$. 
\epF{Lemma}

\ble
[similar to Lemma~\ref{31e}]
\lam{51e}
If\/ $X\sq\Ua$ is a non-empty\/ $\is11$ set then condition\/ $X\ti X$ 
\dd{(\dpx X)}forces that\/ $\doxl\napr\doxr$.\qed
\ele
\vyk{
\bpf
Suppose to the contrary that a subcondition 
$Y\ti Z$ \dd{(\dpx X)}forces $\doxl\apr\doxr$; 
note that 
$Y,Z\sq X$ are non-empty $\is11$ sets with 
$\ek{Y}{\rE_X}=\ek{Z}{\rE_X}$.
Then any \dd{(\dpx X)}generic pair $\ang{y,z}\in Y\ti Z$ 
satisfies $y\apr z$. 
It follows that any two \dd\gh generic reals $y,y'\in Y$ with 
$y\rE_{X} y'$ satisfy $y\apr y'$, too, because one can define 
then a single $z\in Z$ such that both pairs $\ang{y,z}$ and 
$\ang{y',z}$ are \dd{(\dpx X)}generic.

Now we claim that the $\is11$ set 
$W=\ens{\ang{y,y'}\in Y\ti Y}{y\rE_X y'\land y'\napr y}$ 
is empty. 
Indeed, otherwise forcing $\gh_2(W)$ of all non-empty $\is11$ 
sets $P\sq W$ forces pairs of \dd\gh generic (separately) 
reals $y,y'\in Y$ which belong to $W$, hence, satisfy 
$y'\rE_X y$ and $y'\napr y$, which contradicts the above.
Thus if $y,y'\in X$ then $y\rE_X y'\imp y'\apr y$, 
contrary to Lemma~\ref{51-4}.
\epF{Lemma}
}

\ble
\lam{51m}
Let\/ $X\sq \Ua$ be a non-empty\/ $\is11$ set. 
Then\/ $X\ti X$ 
does {\bfit not}\/ \dd{(\dpx X)}force that\/ $\doxl,\doxr$ 
are\/ \dd\cle comparable.
\ele
\bpf
Suppose to the contrary that  
$X\ti X$ \dd{(\dpx X)}forces the comparability, hence,  
forces that $\doxl\cl\doxr\lor\doxl\cl\doxr$ 
by Lemma~\ref{51e}.
There is a 
condition $Y\ti Z$ in $\dpx X$ with $Y\cup Z\sq X$, 
which forces $\doxl\cl\doxr$ or forces $\doxr\cl\doxl$. 
In the second case $Z\ti Y$ forces $\doxl\cl\doxr$ by 
the symmetry.
Thus we can assume that $Y\ti Z$ \dd{(\dpx X)}forces 
$\doxl\cl\doxr$.\vom

1. 
{\it We claim that if\/ $y\in Y,$ $z\in Z$, and\/ $y\rE_X z$ 
then\/ $y\cl z$ strictly\/}.
Indeed if this is not the case then the set 
$$
W=\ens{\ang{y,z}\in Y\ti Z}{y\rE_X z\land y\not\cl z}
$$
is non-empty. 
Let $Y'=\dom W$.
As $Y'\sq X,$ 
by Lemma~\ref{51e} there is a condition $A\ti B\in\dpx X$, 
with $A\cup B\sq Y'$, which forces $\doxl\cl\doxr$; 
for if it forces $\doxr\cl\doxl$ then just consider $B\ti A$ 
instead of $A\ti B$.
By Lemma~\ref{r6} we can assume that  $A\ti B$ is saturated: 
$\ek{A}{\rE_X}=\ek{B}{\rE_X}$.
Consider the forcing notion $\cP$ of all non-$\pu$ $\is11$ sets
of the form $P\ti B'$, where 
$$
P\sq W,\;\;\dom P\sq A,\;\; B'\sq B, 
\qand 
\ek{B'}{\rE_X}=\ek{\dom P}{\rE_X}=\ek{\ran P}{\rE_X}\,.
$$
For instance if $B'=B$ and $P=\ens{\ang{x,y}\in W}{x\in A}$ then 
$P\ti B'\in\cP$. 
Then $\cP$ forces a pair $\ang{\doxl,\doxr}\in W$ and a separate 
real $\dox\in B$ such that both pairs $\ang{\doxl,\dox}$ and 
 $\ang{\doxr,\dox}$ are \dd{(\dpx X)}generic. 
It follows that $\cP$ forces 
both $\doxl\cl\dox$ (as this pair belongs to $A\ti B$) 
and $\dox\cl\doxr$ (as this pair belongs to $Y\ti Z$), 
hence, forces  $\doxl\cl\doxr$. 
On the other hand $\cP$ forces $\doxl\not\cl\doxr$ 
(as this pair belongs to $W$), which is a contradiction.\vom

2. 
It follows by Lemma~\ref{i1} that 
there is a\/ $\id11$ set\/ $\wY\sq\bn$ such that $Y\sq\wY,$  
still\/ $Z\cap\wY=\pu$, and\/ $\wY$ is downwards\/ \dd\cle closed 
in each\/ \dd{\rE_X}class.\vom
 
3. 
We claim that, moreover, 
{\it if\/ $y\in \wY\cap X$, $z\in X\bez\wY$, and\/ $y\rE_X z$, 
then\/ $y\cl z$}. 
Indeed otherwise, the following $\is11$ set 
$$
H_0=\ens{z\in X\bez\wY}
{\sus y\in\wY\cap X\:(y\rE_X z\land y\not\cl z)}
\sq X
$$          
is non-$\pu$. 
As above, 
there is a saturated condition $H\ti H'$ in $\dpx X$, 
with $H\cup H'\sq H_0$, 
which forces $\doxl\cl \doxr$. 
Similarly to part 1 of the proof,  
$z\cl z'$ holds whenever $\ang{z,z'}\in H\ti H'$  and $z\rE_X z'$. 
By construction the $\is11$ set 
$$
\wY_1 =\ens{y\in\wY\cap X}
{\sus z'\in H'\,({y\rE_X z'} \land y\not\cl z')} 
$$
satisfies $\ek{\wY_1}{\rE_X}=\ek{H}{\rE_X}=\ek{H'}{\rE_X}$, 
hence $\wY_1\ti H$ is a condition in $\dpx X$.
Let $\ang{y_1,z}\in \wY_1\ti H$ be any \dd{(\dpx X)}generic pair. 
Then ${y_1\rE_X z}$ by Corollary~\ref{Esmu4}, 
and, by the choice of $X$ and  Lemma~\ref{51e}, 
we have $y_1\cl z$ or $z\cl y_1$. 
However by construction $y_1\in\wY$, $z\nin \wY$, and $\wY$ 
is downwards closed in each \dd{\rE_X}class. 
Thus in fact $y_1\cl z$.
Therefore, for all $z'\in H'$, if $y_1\rE_X z'$ then 
$y_1\cl z\cl z'$, which contradicts to $y_1\in\wY_1$.\vom 

4. 
We conclude, by Lemma~\ref{Lef}, that 
{\it there is a single function\/ $F\in\cF_X$ 
such that if\/ $y\in \wY\cap X$, $z\in X\bez\wY$, and\/ 
$F(y)=F(z)$, then\/ $y\cl z$}.\vom
 
5. 
Following the proof of Corollary~\ref{i2}, 
we claim that {\it the derived  function}
$$
G(x)=\left\{
\bay{rcl}
F(x)\we 0\,,&\text{whenewer}& x\in\wY\\[1ex]
F(x)\we 1\,,&\text{whenewer}& x\in\bn\bez\wY 
\eay
\right.
$$
{\it belongs to\/ $\cfx X.$} 
First of all, still $G\in\cF$ since $C$ is downwards 
\dd\cle closed in each \dd{\dpx X}class. 
Now suppose that $z,y\in X$ and $G(y)\lex G(z)$.
Then either $F(y)\lex F(z)$, or $F(z)=F(y)$ and $y\in\wY$ but 
$z\nin\wY$.
In the ``either'' case immediately $y\cl z$ since $F\in\cfx X$ 
\snos
{The family $\cF$ would not work in the passage; here  
we have to use $\cfx X$ instead.}
. 
In the ``or'' case we have $y\cl z$ by the choice of $F$ 
and the definition of $G$.
Thus $G\in\cfx X$.\vom 

6. 
Pick any pair of reals $y\in Y$ and $z\in Z$ with $y\rE_X z$. 
Then $G(x)=G(y)$ since $G\in\cfx X$. 
But we have $y\in \wY$ and $z\nin\wY$ since 
$Y\sq\wY$ and $Z\cap\wY=\pu$ by construction, and in this case 
surely $G(y)\ne G(z)$ by the definition of $G$.
This contradiction completes the proof of Lemma~\ref{51m}.
%
\epf
}

Lemma~\ref{51+} plus Lemma~\ref{29} imply claim \ref{mtd} 
of Theorem~\ref{mt}.\vtm


\qeDD{Theorem~\ref{mt}}


{\small

\begin{thebibliography}{10}

\bibitem{hkl}
L.~A. Harrington, A.~S. Kechris, and A.~Louveau.
\newblock A {G}limm-{E}ffros dichotomy for {B}orel equivalence relations.
\newblock {\em J. Amer. Math. Soc.}, 3(4):903--928, 1990.

\bibitem{hms}
L.~A. {Harrington}, D.~{Marker}, and S.~{Shelah}.
\newblock {Borel orderings.}
\newblock {\em {Trans. Am. Math. Soc.}}, 310(1):293--302, 1988.

\bibitem{k:blin}
Vladimir Kanovei.
\newblock When a partial {B}orel order is linearizable.
\newblock {\em Fund. Math.}, 155(3):301--309, 1998.

\bibitem{kanB}
Vladimir Kanovei.
\newblock {\em {Borel equivalence relations. Structure and classification.}}
\newblock Providence, RI: American Mathematical Society (AMS), 2008.

\bibitem{ksz}
Vladimir Kanovei, Martin Sabok, and Jind\v{r}ich Zapletal.
\newblock {\em {Canonical Ramsey theory on Polish spaces.}}
\newblock Cambridge: Cambridge University Press, 2013.
\end{thebibliography} 
}


\end{document}